\documentclass[12pt]{article}

\usepackage{amsmath}
\usepackage{amsthm}
\usepackage{makeidx}
\usepackage{latexsym}
\usepackage{amssymb}
\usepackage{enumerate}

\newtheorem{thm}{Theorem}[section]

\newtheorem{lem}[thm]{Lemma}
\newtheorem{prop}[thm]{Proposition}

\theoremstyle{definition}

\newtheorem{rem}[thm]{Remark}

\numberwithin{equation}{section}

\begin{document}

\baselineskip=17pt

\title{Linear recurrence sequences and periodicity of multidimensional continued fractions}

\author{Nadir Murru\\
Department of Mathematics\\ 
University of Turin\\
10123 Turin, Italy\\
E-mail: nadir.murru@unito.it
}

\date{}

\maketitle

\renewcommand{\thefootnote}{}



\renewcommand{\thefootnote}{\arabic{footnote}}
\setcounter{footnote}{0}

\begin{abstract}
Multidimensional continued fractions generalize classical continued fractions with the aim of providing periodic representations of algebraic irrationalities by means of integer sequences. However, there does not exist any algorithm that provides a periodic multidimensional continued fraction when algebraic irrationalities are given as inputs. In this paper, we provide a characterization for periodicity of Jacobi--Perron algorithm by means of linear recurrence sequences. In particular, we prove that partial quotients of a multidimensional continued fraction are periodic if and only if numerators and denominators of convergents are linear recurrence sequences, generalizing similar results that hold for classical continued fractions. 
\end{abstract}

\section{Introduction}
Ternary continued fractions generalize classical continued fractions and can be used to represent pair of real numbers. In particular, a pair of real numbers $(\alpha_0, \beta_0)$ is represented by two sequences $(a_n)_{n=0}^\infty$ and $(b_n)_{n=0}^\infty$ as follows:
\begin{equation} \label{eq:tcf}
\alpha_0=a_0+\cfrac{b_1+\cfrac{1}{a_2+\cfrac{b_3+\cfrac{1}{\ddots}}{a_3+\cfrac{\ddots}{\ddots}}}}{a_1+\cfrac{b_2+\cfrac{1}{a_3+\cfrac{\ddots}{\ddots}}}{a_2+\cfrac{b_3+\cfrac{1}{\ddots}}{a_3+\cfrac{\ddots}{\ddots}}}} \quad \text{and} \quad \beta_0=b_0+\cfrac{1}{a_1+\cfrac{b_2+\cfrac{1}{a_3+\cfrac{\ddots}{\ddots}}}{a_2+\cfrac{b_3+\cfrac{1}{\ddots}}{a_3+\cfrac{\ddots}{\ddots}}}}.
\end{equation}
\emph{Partial quotients} $a_i$ and $b_i$ can be evaluated by means of the Jacobi algorithm:
$$\begin{cases} a_n=[\alpha_n] \cr b_n=[\beta_n] \cr \alpha_{n+1}=\cfrac{1}{\beta_n-[\beta_n]} \cr \beta_{n+1}=\cfrac{\alpha_n-a_n}{\beta_n-b_n}  \end{cases},\quad n=0,1,2,...$$
Ternary continued fractions and previous algorithm were introduced by Jacobi \cite{Jac} in order to answer to a question posed by Hermite. Indeed, Hermite asked for finding a way to give a periodic representation for cubic irrationalities, simlarly to the case of quadratic irrationalities and continued fractions \cite{Her}. This is still a fascinating open problem in number theory, since periodicity of the Jacobi algorithm has been never proved when inputs are cubic irrationals. For this reason, ternary continued fractions have been widely studied during the years. Convergence and properties of ternary continued fractions were studied by Daus and Lehmer in several classical works as \cite{Daus1}, \cite{Daus2}, \cite{Lehmer}. More recently, ternary continued fractions and Jacobi algorithm have been studied for approaching the Hermite problem. In particular several authors tried to modify the Jacobi algorithm in order to obtain a periodic algorithm when cubic irrationalities are given as inputs, \cite{Wil}, \cite{Tam}, \cite{Sch}.

Perron \cite{Perron} generalized the Jacobi approach to higher dimensions, representing a $m$--tuple of real numbers $\alpha^{(1)}_0,...,\alpha^{(m)}_0$ by means of a multidimensional continued fraction defined by $m$ sequences of integers $(a_n^{(1)})_{n=0}^\infty$, ..., $(a_n^{(m)})_{n=0}^\infty$. Partial quotients $a_n^{(1)}$, ..., $a_n^{(m)}$ can be evaluated by means of
\begin{equation} \label{eq:jpa} \begin{cases} a_n^{(i)}=[\alpha_n^{(i)}], \quad i=1,...,m \cr \alpha_{n+1}^{(1)}= \cfrac{1}{\alpha_{n}^{(m)}-a_n^{(m)}} \cr \alpha_{n+1}^{(i)}=\cfrac{\alpha_{n}^{(i-1)}-a_n^{(i-1)}}{\alpha_{n}^{(m)}-a_n^{(m)}} \quad i=2,...,m \end{cases},\quad n=0,1,2,... \end{equation}

Multidimensional continued fractions have been introduced with the aim of finding a periodic representation for all algebraic irrationalities. However previous algorithm (and modifications) has never been proved periodic when algebraic irrationals of order $m$ are given as inputs. Further results and modifications of the Jacobi--Perron algorithm can be found in \cite{Ber}, \cite{Raju}, \cite{Hen}. 

In this paper we show properties that could open new ways for approaching the Hermite problem and proving the periodicity of the Jacobi--Perron algorithm. These properties connect periodicity of partial quotients with linear recurrence sequences. In particular we will prove that sequences $(a_n^{(i)})_{n=0}^\infty$, for $i=1,...,m$, are periodic if and only if sequences of numerators and denominators of convergents of multidimensional continued fractions are recurrence sequences with constant coefficients, generalizing a similar result provided by Lenstra and Shallit for classical continued fractions \cite{Len}. In section \ref{sec:mul}, we give some preliminary properties and notation about multidimensional continued fractions. Section \ref{sec:main} shows the main result and Section \ref{sec:fut} describes further developments of the present work.

\section{Multidimensional continued fractions} \label{sec:mul}
In the following, we write the multidimensional continued fraction derived by equations \eqref{eq:jpa} as follows:
\begin{equation} \label{eq:mcf} (\alpha_0^{(1)}, ..., \alpha_0^{(m)})=[\{a_0^{(1)},a_1^{(1)},a_2^{(1)},...\},...,\{a_0^{(m)},a_1^{(m)},a_2^{(m)},...\}]. \end{equation}
This multidimensional continued fraction generalizes the structure of ternary continued fraction \eqref{eq:tcf} by means of the following relations
$$\begin{cases}  \alpha_n^{(i-1)}=a_n^{(i-1)}+\cfrac{\alpha^{(i)}_{n+1}}{\alpha_{n+1}^{(1)}}, \quad i=2,...,m \cr \alpha_n^{(m)}=a_n^{(m)}+\cfrac{1}{\alpha_{n+1}^{(1)}} \end{cases},\quad n=0,1,2,...,$$
which can be derived by equations \eqref{eq:jpa}. Let us observe that a multidimensional continued fraction is finite if and only if $\alpha_0^{(1)},...,\alpha_0^{(m)}$ are rational numbers, similarly to classical continued fractions. In the following, we focus on $\alpha_0^{(1)},...,\alpha_0^{(m)}$ irrational numbers.

We can introduce the notion of \emph{convergents}, like for classical continued fractions. In this case, we have $m$ sequences of convergents $\left( \cfrac{A_n^{(i)}}{A_n^{(m+1)}}\right)_{n=0}^\infty$, for $i=1,...,m$, defined by
\begin{equation} \label{eq:conv0} \left( \cfrac{A_n^{(1)}}{A_n^{(m+1)}}, ..., \cfrac{A_n^{(m)}}{A_n^{(m+1)}} \right)=[\{a_0^{(1)},...,a_n^{(1)}\},...,\{a_0^{(m)},...,a_n^{(m)}\}], \quad \forall n\geq0 \end{equation}
and we have
$$\lim_{n\rightarrow\infty}\cfrac{A_n^{(i)}}{A_n^{(m+1)}}=\alpha_0^{(i)},\quad i=1,...,m.$$
In \cite{Ber1}, the reader can find useful properties regarding convergents of multidimensional continued fractions (in this paper we have used a little different notation). In particular, we highlight the following properties.
\begin{prop}
Numerators and denominators of convergents of multidimensional continued fraction \eqref{eq:mcf} satisfy the following recurrences of order $m+1$ with non--constant coefficients:
\begin{equation} \label{eq:conv} A_n^{(i)}=\sum_{j=1}^m{a_n^{(j)}A_{n-j}^{(i)}}+A_{n-m-1}^{(i)},\quad i=1,...,m+1, \quad \forall n\geq 1 \end{equation} 
with initial conditions
$$\begin{cases} A_0^{(i)}=a_0^{(i)},  \quad  A_0^{(m+1)}=1, \cr A_{-i}^{(i)}=1, \quad  A_{-i}^{(m+1)}=0, \cr A_{-j}^{(i)}=0, \quad j=1,...,m, \quad j\not=i \end{cases}\quad i=1,...,m.$$ 
\end{prop}
\begin{prop} \label{prop:mat}
Given the multidimensional continued fraction \eqref{eq:mcf}, we have
$$ \prod_{j=0}^n \begin{pmatrix} a_j^{(1)} & 1 & 0 & 0 & ... & 0 \cr a_j^{(2)} & 0 & 1 & 0 &... & 0 \cr \vdots & \vdots & \ddots & \ddots & \ddots & \vdots \cr a_j^{(m-1)} & 0 & 0 & ... & 1 & 0 \cr a_j^{(m)} & 0 & 0 & 0 & ... & 1 \cr 1 & 0 & 0 & 0 & ... & 0 \end{pmatrix}=\begin{pmatrix} A_n^{(1)} & A_{n-1}^{(1)} & ... & A_{n-m}^{(1)} \cr A_n^{(2)} & A_{n-1}^{(2)} & ... & A_{n-m}^{(2)} \cr \vdots & \vdots & ... & \vdots \cr A_n^{(m+1)} & A_{n-1}^{(m+1)} & ... & A_{n-m}^{(m+1)}  \end{pmatrix}$$
\[ \Updownarrow \] \[\left( \cfrac{A_n^{(1)}}{A_n^{(m+1)}}, ..., \cfrac{A_n^{(m)}}{A_n^{(m+1)}} \right)=[\{a_0^{(1)},...,a_n^{(1)}\},...,\{a_0^{(m)},...,a_n^{(m)}\}], \]
for $n=0,1,2,...$.
\end{prop}

\begin{rem}
Matrices in the previous proposition are $(m+1)\times(m+1)$ matrices.
\end{rem}

\begin{rem}
In the following, when there is no possibility of confusion, we only write $(a_n)$ instead of $(a_n)_{n=0}^\infty$.
\end{rem}

\section{Main result} \label{sec:main}
In this section, we prove that periodicity of Jacobi--Perron algorithm is strictly related to linear recurrence sequences. We see that a multidimensional continued fraction is periodic if and only if numerators and denominators of convergents are linear recurrence sequences. In other words, we prove that sequences $(a_n^{(i)})$, for $i=1,...,m$, are periodic if and only if sequences $(A_n^{(i)})$, for $i=1,...,m+1$, satisfy a recurrence relation with constant coefficients. Proofs follow ideas used in \cite{Len}, but here the situation is complicated by the presence of several sequences of partial quotients.
\begin{thm}
If sequences of partial quotients of \eqref{eq:mcf} are periodic, then sequences of numerators and denominators of its convergents are linear recurrence sequences.
\end{thm}
\begin{proof}
Let us consider the periodic multidimensional continued fraction
$$[\{a_0^{(1)},...,a_{p_1-1}^{(1)},\overline{b_0^{(1)},...,b_{q_1-1}^{(1)}}\},...,\{a_0^{(m)},...,a_{p_m-1}^{(m)},\overline{b_0^{(m)},...,b_{q_m-1}^{(m)}}\}].$$
We introduce the following quantities: 
$$u=lcm(q_1,...,q_m),\quad v_i=max(p_1,...,p_m)-p_i, \quad i=1,...,m.$$
We define $u$ matrices 
$$M_i=\prod_{0\leq j<u} \begin{pmatrix} b_{i+j+v_1}^{(1)} & 1 & 0 & 0 & ... & 0 \cr b_{i+j+v_2}^{(2)} & 0 & 1 & 0 &... & 0 \cr \vdots & \vdots & \ddots & \ddots & \ddots & \vdots \cr b_{i+j+v_{m-1}}^{(m-1)} & 0 & 0 & ... & 1 & 0 \cr b_{i+j+v_m}^{(m)} & 0 & 0 & 0 & ... & 1 \cr 1 & 0 & 0 & 0 & ... & 0 \end{pmatrix},\quad \forall i\in \mathbb Z_u.$$
We can easily observe that $\det(M_i)=(-1)^{(m-1)u}$, furthermore $tr(M_i^k)=tr(M_j^k)$ for all $i,j\in\mathbb Z_u$, $\forall k\in\mathbb N$.
Indeed, given any $i,j\in \mathbb Z_u$, it is always possible to find two matrices $P$, $Q$ such that $M_i=PQ$ and $M_j=QP$, and consequently $tr(M_i)=tr(M_j)$. Similar observation holds for powers of matrices $M_i$. Therefore, all matrices $M_i$ have same invariants and characteristic polynomial
$$x^m-\sum_{i=1}^{m-1}h_ix^i-(-1)^{(m-1)u},$$
where $h_i$ are invariants of these matrices.

By proposition \ref{prop:mat} follows that
$$ \begin{pmatrix} A_{n+u}^{(1)} & A_{n+u-1}^{(1)} & ... & A_{n+u-m}^{(1)} \cr A_{n+u}^{(2)} & A_{n+u-1}^{(2)} & ... & A_{n+u-m}^{(2)} \cr \vdots & \vdots & ... & \vdots \cr A_{n+u}^{(m+1)} & A_{n+u-1}^{(m+1)} & ... & A_{n+u-m}^{(m+1)}  \end{pmatrix}= \begin{pmatrix} A_n^{(1)} & A_{n-1}^{(1)} & ... & A_{n-m}^{(1)} \cr A_n^{(2)} & A_{n-1}^{(2)} & ... & A_{n-m}^{(2)} \cr \vdots & \vdots & ... & \vdots \cr A_n^{(m+1)} & A_{n-1}^{(m+1)} & ... & A_{n-m}^{(m+1)}  \end{pmatrix}M_{n-u}$$
for all $n\geq \max(p_1,...,p_m)$. Now, we can write
$$M_{n-u}^m-\sum_{i=1}^{m-1}h_iM_{n-u}^i-(-1)^{(m-1)u}I=O$$
where $I$ is the identity matrix and $O$ is the zero matrix. Hence, multiplying (on the left) both sides of previous equation by 
$$\begin{pmatrix} A_n^{(1)} & A_{n-1}^{(1)} & ... & A_{n-m}^{(1)} \cr A_n^{(2)} & A_{n-1}^{(2)} & ... & A_{n-m}^{(2)} \cr \vdots & \vdots & ... & \vdots \cr A_n^{(m+1)} & A_{n-1}^{(m+1)} & ... & A_{n-m}^{(m+1)}  \end{pmatrix}$$
we obtain
\begin{align*}
\begin{pmatrix} A_{n+mu}^{(1)} & A_{n+mu-1}^{(1)} & ... & A_{n+mu-m}^{(1)} \cr A_{n+mu}^{(2)} & A_{n+mu-1}^{(2)} & ... & A_{n+mu-m}^{(2)} \cr \vdots & \vdots & ... & \vdots \cr A_{n+mu}^{(m+1)} & A_{n+mu-1}^{(m+1)} & ... & A_{n+mu-m}^{(m+1)}  \end{pmatrix}-\sum_{i=1}^{m-1}h_i\begin{pmatrix} A_{n+iu}^{(1)} & A_{n+iu-1}^{(1)} & ... & A_{n+iu-m}^{(1)} \cr A_{n+iu}^{(2)} & A_{n+iu-1}^{(2)} & ... & A_{n+iu-m}^{(2)} \cr \vdots & \vdots & ... & \vdots \cr A_{n+iu}^{(m+1)} & A_{n+iu-1}^{(m+1)} & ... & A_{n+iu-m}^{(m+1)}  \end{pmatrix}+\\
-(-1)^{(m-1)u}\begin{pmatrix} A_{n}^{(1)} & A_{n-1}^{(1)} & ... & A_{n-m}^{(1)} \cr A_{n}^{(2)} & A_{n-1}^{(2)} & ... & A_{n-m}^{(2)} \cr \vdots & \vdots & ... & \vdots \cr A_{n}^{(m+1)} & A_{n-1}^{(m+1)} & ... & A_{n-m}^{(m+1)}  \end{pmatrix}=O.
\end{align*}
Thus, we have, e.g.,
$$A^{(1)}_{n+mu}-\sum_{i=1}^{m-1}h_iA^{(1)}_{n+iu}-(-1)^{m-1}A^{(1)}_n=0$$
and similar relations for sequences $(A_n^{(i)})$, for $i=2,...,m+1$, i.e., sequences $(A_n^{(i)})$ are all linear recurrence sequences.
\end{proof}

We need some lemmas for proving the vice versa. Next lemma is known as the Hadamard quotient theorem and it is due to Van der Poorten \cite{Van}.

\begin{lem} \label{lemma:hadamard}
Let $\sum_{i=0}^\infty a_ix^i$ and $\sum_{i=0}^\infty b_ix^i$ be formal series representing rational functions in $\mathbb C[x]$, if $\cfrac{a_i}{b_i}$ is an integer for all $i\geq0$, then $\sum_{i=0}^\infty \cfrac{a_i}{b_i}x^i$ is a formal series representing a rational function.
\end{lem}

\begin{rem}
We recall that the generating function of a linear recurrence sequence is a rational function.
\end{rem}

We also need some results about sums and products of linear recurrence sequences. The reader can find an interesting overview on this topic in \cite{Cer}. In the next lemma, we summarize some results.

\begin{lem} \label{lemma:sum-prod}
If $(a_n)$ and $(b_n)$ are linear recurrence sequences with characteristic polynomials $p(x)$ and $q(x)$, then $(a_n+b_n)$ and $(a_n\cdot b_n)$ are also linear recurrence sequences whose characteristic polynomials are $p(x)\cdot q(x)$ and $p(x)\odot q(x)$, respectively, where $p(x)\odot q(x)$ is the characteristic polynomial of the matrix obtained by the Kronecker product between companion matrices of $p(x)$ and $q(x)$.
\end{lem}

The following lemma has been also used in \cite{Len}.

\begin{lem} \label{lemma:ls}
Let $(a_n)$ and $(b_n)$ be linear recurrence sequences such that characteristic polynomial of $(b_n)$ divides characteristic polynomial (of degree $e$) of $(a_n)$. Then, there exists constants $k$ and $n_0$ such that
$$\max(\lvert a_n \rvert,\lvert a_{n-1}\rvert,...,\lvert a_{n-e+2}\rvert,\lvert a_{n-e+1} \rvert)>k\lvert b_n \rvert,$$
for any $n\geq n_0$.
\end{lem}

Finally, we prove the following lemma.

\begin{lem} \label{lemma:bound}
With the above notation, we have
$$A^{(i)}_n\geq \prod_{j=1}^n a_j^{(1)},\quad i=1,...,m+1$$
for $n\geq1$.
\end{lem}
\begin{proof}
The thesis can easily proved by induction. Let us observe that, by equations \eqref{eq:jpa}, we have $a_n^{(i)}>0$, for $i=1,...,m$ and $n\geq 1$. The basis is straightforward:
$$A^{(m+1)}_1=a_1^{(1)},\quad A^{(m+1)}_2=a_2^{(1)}a_1^{(1)}+a_2^{(2)}.$$
For the inductive step, we have
$$A_n^{(m+1)}=\sum_{j=1}^m{a_n^{(j)}A_{n-j}^{(m+1)}}+A_{n-m-1}^{(m+1)}\geq a_n^{(1)}\prod_{j=1}^{n-1}a_j^{(1)}+\sum_{j=2}^m{a_n^{(j)}A_{n-j}^{(m+1)}}+A_{n-m-1}^{(m+1)}.$$
Since $A^{(m+1)}_i$, for $i=1,2,...$, are concordant, the thesis follows. Similar considerations hold for sequences $(A_n^{(i)})$, for $i=1,...,m$.
\end{proof}

\begin{thm}
Let $\left(\cfrac{A_n^{(i)}}{A_n^{(m+1)}} \right)$ be the sequences of convergents, for $i=1,...,m$, of the multidimensional continued fraction \eqref{eq:mcf}. If $(A_n^{(i)})$, for $i=1,...,m+1$, are linear recurrence sequences, then $(a_n^{(i)})$, for $i=1,...,m$, are periodic sequences.
\end{thm}
\begin{proof}
For $n$ sufficiently large $A_n^{(i)}\not=0$, for $i=1,...,m+1$, and by equations \eqref{eq:conv}, with a little bit of calculation, it is possible to obtain
$$a_n^{(m)}=\cfrac{F_n}{G_n},$$
where $F_n$ and $G_n$ are sums and products of some $A_{n-j}^{(i)}$, for certain values of $i$ and $j$. Since $(A_n^{(i)})$ are all linear recurrence sequences, by lemma \ref{lemma:sum-prod}, $(F_n)$ and $(G_n)$ are linear recurrence sequences. Moreover, $(F_n)$ and $(G_n)$ satisfy hypothesis of lemma \ref{lemma:hadamard} and consequently $(a_n^{(m)})$ is a linear recurrence sequence. Similarly, we can prove that $(a_n^{(i)})$, for $i=1,...,m-1$ are linear recurrence sequences.

Now, we prove that sequence $(a_n^{(1)})$ is bounded. Let us suppose that characteristic polynomial of the linear recurrence sequence $(a_n^{(1)})$ has degree $e$. By the Binet formula we can write
\begin{equation}\label{eq:binet} a_n^{(1)}=\sum_{i=1}^ep_i(n)\alpha_i^n, \quad n=0,1,2,... \end{equation}
where $\alpha_i$'s are roots of the characteristic polynomial and $p_i(n)$ polynomials of degree $t_i$ in $n$, for $i=1,...,e$. By lemma \ref{lemma:bound} we have $A^{(m+1)}_{qe}\geq \prod_{1\leq j \leq qe} a_j^{(1)}$, for $q\geq1$. Moreover, for a given index $i$, sequences $(a_n^{(1)})$ and $(n^{t_i}\alpha_i^n)$ satisfy hypothesis of lemma \ref{lemma:ls}. 

Hence by lemma \ref{lemma:ls}, fixed $1\leq h \leq m$, we have $m$ sequences $a^{(1)}_{(h-1)e+1},...,a^{(1)}_{he}$, such that the maximum element is greater than $k(he)^{t_i}\lvert \alpha_i^e \rvert^h$. Said $k'$ a constant greater than the product of all other elements, we have
$$A^{(m+1)}_{qe}>k'k^m(m!)^{t_i}e^{mt_i}\lvert \alpha_i^e \rvert^{m(m+1)/2}.$$
A similar inequality was find in \cite{Len} for denominators of convergents of classical continued fractions.

Since $(A^{(m+1)}_n)$ is a linear recurrence sequence, it can not grow too rapidly. Indeed, it is well--known that $\log A^{(m+1)}_n =O(n)$ and consequently $\lvert \alpha_i \rvert <1$ or $\lvert \alpha_i\rvert=1$ and $t_i=0$. Thus, by equation \eqref{eq:binet}, we have that $(a_n^{(1)})$ is bounded. Since $(a^{(1)}_n)$ is a bounded linear recurrence sequence of integers, we have proved that $(a^{(1)}_n)$ is periodic. Finally, by equations \eqref{eq:jpa} we have that if $(a^{(1)}_n)$ is periodic then sequences $(a^{(i)}_n)$, for $i=2,...,m$, are periodic.
\end{proof}

\section{Future developments} \label{sec:fut}
We have found a characterization of periodicity for the Jacobi--Perron algorithm by means of linear recurrence sequences. This result could be exploited for approaching the Hermite problem by new ways. 

In \cite{Mur}, the author found a periodic representation for any cubic irrational using a ternary continued fraction with rational partial quotients, whose convergents have numerators and denominators that are linear recurrence sequences. However, this representation depends on the knowledge of the minimal polynomial of the involved cubic irrational and consequently none algorithm can be found for generating such a representation. In particular, given the cubic irrational $\alpha$, whose minimal polynomial is $x^3-px^2-qx-r$, we have
\begin{multline}\label{eq:murru} [\{z,\cfrac{2z+p^2+q}{pq+r},\overline{\cfrac{(pq+r)\text{tr}(N)}{\det (N)},\text{tr}(N),\cfrac{\text{tr}(N)}{pq+r}}\}, \\ \{p,-\cfrac{z^2+qz+p^2z-pr}{pq+r},\overline{-\cfrac{I_1(N)}{\det (N)},-\cfrac{(pq+r)I_1(N)}{\det(N)},-\cfrac{I_1(N)}{pq+r}}\}] = 
(\cfrac{r}{\alpha},\alpha)  \end{multline}
where $I_1$ is the second invariant of a $3\times3$ matrix. Moreover, 
$$N=\begin{pmatrix} z & r & pr \cr 0 & q+z & pq+r \cr 1 & p & p^2+q+z \end{pmatrix} $$
for a given integer $z$.

Let us define
\begin{equation} \label{eq:noper} [\{a_0,a_1,...\},\{b_0,b_1,...\}]=(\cfrac{r}{\alpha},\alpha)  \end{equation}
the ternary continued fraction whose partial quotients are derived by the Jacobi algorithm. Since ternary continued fractions \eqref{eq:murru} and \eqref{eq:noper} are equal, it could be interesting to compare convergents of \eqref{eq:murru} with convergents of \eqref{eq:noper} with the aim of proving that numerators and denominators of convergents of \eqref{eq:noper} are linear recurrence sequences too.

\section{Ackonwledgments}
I would like to thank Prof. Umberto Cerruti for introducing me to these beautiful problems.

\end{document}